\documentclass[11pt]{article}
\usepackage{times}
\usepackage{amsmath,amstext,amssymb,amsthm, amsfonts}
\usepackage[dvips]{graphics,color}

\newcommand{\Pf}{\paragraph{{\bf Proof.}}}       
\newcommand{\blot}{\hfill{\vrule height .9ex width .8ex depth -.1ex }}
\newcommand{\EndPf}{\hfill $\blot$ \medskip}     
\newcommand{\bi}{\begin{itemize}}  
\newcommand{\ei}{\end{itemize}}     
\newcommand{\be}{\begin{enumerate}}  
\newcommand{\ee}{\end{enumerate}}
\newcommand{\bc}{\begin{center}}  
\newcommand{\ec}{\end{center}}     
\def\lm{\lambda}
\def\tlr{\tilde{r}}
\def\tlR{\tilde{R}}
\def\tlq{\tilde{q}}
\def\tlQ{\tilde{Q}}
\def\tlV{\tilde{V}}

\def\rsb{\right ]}
\def\lsb{\left [}
\def\lcb{\left \{}
\def\rcb{\right \}}
\def\argmin{\mathop{\hbox{\rm  argmin}}}
\newcommand{\beq}{\begin{equation}}
\newcommand{\eeq}{\end{equation}}
\newcommand{\ls}[1]
   {\dimen0=\fontdimen6\the\font \lineskip=#1\dimen0
   \advance\lineskip.5\fontdimen5\the\font \advance\lineskip-\dimen0
   \lineskiplimit=.9\lineskip \baselineskip=\lineskip
   \advance\baselineskip\dimen0 \normallineskip\lineskip
   \normallineskiplimit\lineskiplimit \normalbaselineskip\baselineskip
   \ignorespaces }

\numberwithin{equation}{section}

\title{The Multi-Armed Bandit, with Constraints}

\begin{document}

\newtheorem{Lemma}{Lemma}[section]
\newtheorem{Theorem}[Lemma]{Theorem}
\newtheorem{Corollary}[Lemma]{Corollary}
\newtheorem{Definition}[Lemma]{Definition}
\newtheorem{prop}[Lemma]{Proposition}
\newtheorem{Ex}[Lemma]{Example}
\newtheorem{Alg}[Lemma]{Algorithm}
\newtheorem{Remark}[Lemma]{Remark}
\newtheorem{Assumption}[Lemma]{Assumption}

\date{}
\maketitle

\vspace{-2cm}

\begin{center}
Eric V. Denardo,\footnote{Center for Systems Sciences, Yale University, PO Box 208267, New Haven, CT 06520, USA.} \and
Eugene A. Feinberg
\footnote{Department of Applied Mathematics and Statistics,
  Stony Brook University,
Stony Brook, NY 11794-3600, USA.} and Uriel G.
Rothblum\footnote{Faculty of Industrial Engineering and
Management, Technion -- Israel Institute of Technology, Haifa
32000, Israel.} \\

\bigskip
\date{\today}

\end{center}

\begin{abstract}
\noindent
The early sections of this paper present an analysis of a Markov decision model that is known as the multi-armed bandit under the assumption that the utility function of the decision maker is either linear or exponential. The analysis includes efficient procedures for computing the expected utility associated with the use of a priority policy and for identifying a priority policy that is optimal. The methodology in these sections is novel, building on the use of elementary row operations. In the later sections of this paper, the analysis is adapted to accommodate constraints that link the bandits.
\end{abstract}

\section{Introduction} \label{s1-111113}
The colorfully-named \emph{multi-armed bandit} \cite{G79} is the following Markov decision problem:  At epochs  $1,  2, \dots$,  a decision maker observes the current state of each of several Markov chains with rewards (bandits) and plays one of them.  The Markov chains that are not played remain in their current states.  The Markov chain that is played evolves for one transition according to its transition probabilities, earning an immediate reward (possibly negative) that can depend upon its current state and on the state to which transition occurs.  Henceforth, to distinguish the states of the individual Markov chains from those of the Markov decision problem, the latter are called multi-states; each multi-state prescribes a state for each of the Markov chains.

A key result for the multi-armed bandit is that attention can be restricted to a simple class of decision procedures that are based on ``labelings."  A {\em labeling}\/ is an assignment of a number to each state of each bandit such that no two states have the same number (label), even if they are in different bandits.  A {\em priority rule}\/ is a policy that is determined by a labeling in this way; given each multi-state, the priority rule plays the Markov chain whose current state has the lowest label.
In a seminal 1974 paper, Gittins and Jones \cite{GJ74} (followed by \cite{G79}) demonstrated the optimality of a priority rule for a model whose objective is to maximize expected discounted income with a per-period discount factor  $c$  having  $0 < c  <  1$. The (optimal) priorities that they identified are based on a family of stopping times, one for each state of each chain.  Given state  $i$  of bandit  $k$,  the decision maker is imagined to play bandit  $k$  for any number
$\tau$  ($\tau\ge 1$)  of consecutive epochs, observing the state to which each transition occurs, and stopping whenever he or she wishes to do so.  The discounted present value of the (random) income stream that is received during epochs  1  through  $\tau$   is denoted $X(\tau)$.  The stopping times $\tau$   for state  $i$  are used to assign that state an {\em index}\/  $I(i)$  by
\beq
\label{eq2.1}		
I(i) = \max_{\tau} \lcb \frac{E[X(\tau]}{1-E[c^{\tau}]} \rcb .
\eeq
It was demonstrated in \cite{GJ74,G79} that, given each multi-state, it is optimal to play any Markov chain (bandit) whose current state has the largest index (lowest label).

	Following \cite{GJ74,G79}, the multi-armed bandit problem has stimulated research in control theory, economics, probability, and operations research.  A sampling of noteworthy papers includes Bergemann and V\"{a}lim\"{a}kim \cite{berg08}, Bertsimas and Ni\~{n}o-Mora \cite{B-NM93}, El Karoui and Karatzas \cite{elka94}, Katehakis and Veinott \cite{kate87}, Schlag \cite{schl98}, Sonin \cite{soni08}, Tsisiklis \cite{tsit94}, Variaya, Walrand and Buyukkoc \cite{vari85}, Weber \cite{webe92}, and Whittle \cite{whit80},. Books on the subject (that list many references) include Berry and Fristedt \cite{BF85}, Gittins \cite{G89}, Gittins, Glazebrook and Weber \cite{GGW11}. The last and most recent of these books provides a status report on the multi-armed bandit that is almost up-to-date.

	An implication of the analysis in \cite{GJ74,G79} is that the largest of all of the indices equals the maximum over all states of the ratio $r(i)/(1-c)$,  where $r(i)$
denotes the expectation of the reward that is earned if state  $i$'s  bandit is played once while state  $i$  is observed and where  $c$  is the discount factor.
In 1994, Tsitsiklis \cite{tsit94} observed that repeated play of a bandit while it is
in the state  $i $ whose ratio is largest leads to a multi-armed bandit with one
fewer state and random transition times.

	In 2007 Denardo, Park and Rothblum \cite{DPR}  considered a generalization of classic multi-armed bandit model with the following new features:
\bi
\item
The utility function of the decision maker can be exponential, expressing sensitivity to risk.
\item
In the case of linear utility functions, the assumption that rewards are discounted is replaced by the introduction of stopping (which captures discounting).
\ei
The analysis of \cite{DPR} focused on pair-wise comparisons  It avoided the use of stopping times, which had been a common feature of the prior analyses of multi-armed bandits. It relied on linear algebra, rather than on probability theory.  It avoided the need to deal, in the more general cases, with ratios that had zeros in their denominators. It included efficient algorithms for computing indices and for identifying an optimal priority rule.

 Constraints that link the bandits (for the extension considered in \cite{DPR}) are dealt with in Sections 7-8 of the current paper. An
 optimal solution to the multi-armed bandit problem with $W$ constraints
 is shown to be an initial randomization over $W+1$ priority rules, each of which
 is the optimal solution to an unconstrained bandit problem whose rewards are determined by a particular set of prices (multipliers) on the constraints. A column generation algorithm
 is described for computing such an optimal solution. In each stage, the coefficients of the
 column that enters the basis are found by the application of the policy evaluation procedure of Section 4.

As concerns contributions to methodology, the analysis in the earlier sections of this paper rests on elementary row operations.  Row operations are used in sections 3-4 to present the first efficient algorithm for computing the utility function gained when beginning at a given multi-state and using any given priority rule (solving the optimality equation is inefficient as the number of multi-states can be enormous). Row operations are also used in Section 5-6 to determine efficiently an optimal priority policy and to provide a proof for its optimality.
The approach in the current paper builds on that of \cite{DPR}, but simplifies the theoretical development
and the computation. In particular, the computation effort that the method requires matches the best existing bound for computing Gittins indices (obtained in \cite{NM07}, see \cite[p.43]{GGW11}  (in fact, the same bound applies to the method developed in \cite{DPR}).

\section{The model}\label{s3-111113}
Let  $K$  be the number of Markov chains (bandits), and let them be numbered  $1$  through $K$.  Markov chain  $k$  has a finite set  $N_k$  of states.  No loss of generality occurs by
assuming, as we do, that the states of distinct Markov chains are disjoint.  Thus, each state $j$  identifies the Markov chain  $b(j)$  of which it is a member, i.e., $j \in N_{b(j)}$.  The set of all states of all bandits is given by  $N = N_1 \cup \dots \cup N_K$.

If bandit  $k$  is played while its state is  $i$,  this bandit experiences transition to state  $j$  with probability  $p(i, j)$,  and it experiences \emph{termination of play} with probability $p(i,0)$ given by
\[
p(i,0) = 1  - \sum_{j \in N_k} p(i,j) \qquad \forall \, i \in N_k, \quad \forall \, k \in \{1,2,\dots, K\} \, . \]   		
If bandit  $k$  is played when its state is  $i$  and if transition is to occur to state  $j$,  payoff  $x(i, j)$  is earned at the start of the period; if termination is to occur instead, payoff  $x(i, 0)$  is earned at the start of the period.  Each of these ``payoffs" can be positive, negative or zero.

Termination stops the play of all  $K$  bandits, not merely of the bandit that is being played.  Termination is modeled as transition to state  $0$.  No action is possible after transition to state  $0$.  For this reason, state  $0$  is excluded from  $N_k$   for each  $k$,  and hence from $N$.

Each multi-state  $s$  is a set that contains, for each  $k$,  exactly one state in $N_k$.  When  $s$  is a multi-state, the symbol  $s_k$   denotes the state of bandit  $k$  that is included in  $s$,  and the symbol  $s_{\setminus k}$ is defined by  $s_{\setminus k} = s \setminus \{s_k\}$.  Thus,  $s_{\setminus k}$    contains all the states in  $s$ other than $s_k$.  Let  $S$  denote the set of all multi-states.  Given any multi-state  $s$,  one of the bandits must be played.  Hence, for this model, a \emph{stationary nonrandomized policy}  $\delta$  is any map that for each multi-state $s \in S$ picks a bandit  $\delta (s) \in \{1,2,\dots, K\}$.  Let $\Delta$ denote the set of all stationary nonrandomized policies.

\subsection{Utility}

The goal is to maximize expected utility.  This will be accomplished with a \emph{linear utility function}  $u(x)  =  x$,  with a \emph{risk-averse exponential utility function} $u(x) = -e^{-\lm x}$  where $\lm$  is a positive constant that is known as the coefficient of risk aversion and with a \emph{risk-seeking exponential utility function} $u(x) = e^{\lm x}$ where $\lm$  is a positive constant.

All three cases are described and analyzed using the \emph{local utility function}  $h(s, k, v)$  whose value equals the expectation of the total utility that is earned in the (artificially-truncated) one-transition model if multi-state  $s$  is observed now, if bandit  $k$  is selected now, and if utility $v(t)$   is earned if transition occurs to multi-state $t$.

\subsection{Linear utility}

In the case of linear utility, the local utility function is
\beq
\label{eq3.1}  h(s, k, v)  = r(s_k)    +  \sum_{j \in N_k} q(s_k,j) v(s_{\setminus k} \cup \{j\}) \, ,	
\eeq
with data  $r(i)$  and  $q(i,j)$   that are specified by
\beq
\label{eq3.2}    	r(i)  =   p(i,0) x(i, 0)  +  \sum_{j \in N_{b(i)}} p(i,j) x(i,j) \qquad \forall \, i \in N \, ,
\eeq
\beq
\label{eq3.3} 		q(i, j)  =  p(i, j) \, ,	 \qquad \forall \, i, j \in N_k \, , \quad \forall \, k \in \{1,2,\dots, K\} \, .
\eeq
Interpret  $r(i)$  as the expectation of the reward that is earned immediately if bandit  $b(i)$  is played while its state is  $i$,  and interpret $q(i,j)$    as the probability that bandit  $b(i)$  will experience transition to state  $j$  given that it is played while its state is  $i$.  As noted earlier, playing bandit  $b(i)$  while its state is $ i$  causes termination (rather than transition to some state  $j$  in   $N_{b(i)}$)  with probability  $p(i, 0)$,  which may be positive.
The above model captures the classic discounted model, which has transition probability  $p(i, j)$ and  discount factor  $c$  satisfying $0  <  c  <  1$, by replacing $p(i,0)$ and
$p(i,j)$ in (\ref{eq3.2})-(\ref{eq3.3}) by $cp(i,0)$ and
$cp(i,j)$, respectively. Incorporating the discount factor into the transition rates yields a fundamental advantage -- it facilitates
an analysis that applies linear algebraic arguments instead of stopping times.

\subsection{Exponential utility}

With the risk-averse exponential utility function $u(x)  = -e^{-\lm x}$,  one has  $u(x + y)  = -e^{-\lm (x+y)} =  e^{-\lm x} u(y)$,  and the local utility function is given by (\ref{eq3.1}) with data  $r(i)$  and  $q(i, j)$  that are specified, for each  $i \in N$  and $j \in N_{b(i)}$,  by
\beq
\label{eq3.4}		r(i) = -p(i,0) \, e^{-\lm x(i,0)} \ \mbox{ and  } \ q(i,j) = p(i,j) \, e^{-\lm x(i,j)} .
\eeq

With the risk-seeking exponential utility function  $u(x)  = e^{\lm x}$,  the local utility function is given by (\ref{eq3.1}) with data
\beq
\label{eq3.5} r(i) = p(i,0) \, e^{\lm x(i,0)} \ \mbox{ and  } \ q(i,j) = p(i,j) \, e^{\lm x(i,j)} .
\eeq

With all three utility functions,  $r(i)$  is called a \emph{reward}, and  $q(i, j)$  is called a \emph{transition rate}.  In the linear-utility model,  $q(i, j)$  is a probability.   In the risk-averse exponential-utility model, $ q(i,j)$ is the product of a probability and a disutility.

Bandit  $k$  has an  $|N_k | \times | N_k|$   matrix  $q^k$   whose  $ij^{\rm th}$ entry equals  $q(i, j)$  for each ordered pair  $(i, j)$  of states in $N_k$.  In each case, the entries in  $q^k$ are nonnegative.  In the linear-utility case, $q^k$ is \emph{substochastic}, which is to say that its entries are nonnegative and the entries in each row sum to  $1$  or less.  In the risk-averse exponential case, each state  $i$  has reward  $r(i) \leq 0$.  In the risk-seeking exponential case, each state  $i$  has  reward  $r(i) \geq 0$.

\subsection{A hypothesis}
 	
 A square matrix  $Q$  is \emph{transient} if and only if each entry in the matrix $Q^t$ approaches $0 $ as   $t \to \infty$.  A hypothesis that is shared by all three utility functions is presented below as:

 \medskip

 \noindent
 {\bf Hypothesis C.} Expressions (\ref{eq3.6}) and at least one of (\ref{eq3.7}), (\ref{eq3.8}) and (\ref{eq3.9}) are satisfied.
\begin{align}
& q^k  \quad \mbox{ is nonnegative and transient} \quad & \mbox{for } \  k  =  1, 2, \dots, K \, . \ \, \
\label{eq3.6} \\
& q^k  \quad \mbox{ is substochastic} & \mbox{for } \  k  =  1, 2, \dots, K \, .  \ \, \
\label{eq3.7} \\		   		
& r(i)  \leq 0 & \mbox{for } \ i=1,2,\dots, |N| \, . \ \label{eq3.8} \\
& r(i)  \geq 0 & \mbox{for } \ i=1,2,\dots, |N| \, . \ \label{eq3.9}
\end{align}

In the case in which (\ref{eq3.6}) and (\ref{eq3.7}) hold is dubbed \emph{Hypothesis RN} (short for risk neutral).  This case includes the classic discounted model, which has transition probability  $p(i, j)$,  discount factor  $c$  that satisfies  $0  <  c  <  1$  and  $q(i, j)  =  c p(i, j)$,  so that each row of $(q_k)^t$ sums to $c$,  which guarantees that $q^k$ is transient.  Hypothesis RN also encompasses linear-utility models in which ratios akin to (\ref{eq2.1}) would have $0$'s  in their denominator.  Hypothesis RN is relaxed in Section 10.

The case in which (\ref{eq3.6}) and (\ref{eq3.8}) hold is dubbed \emph{Hypothesis RA} (short for risk-averse).  In this case, the assumption that $q^k$  is transient excludes a bandit whose repeated play would earn expected utility of $-\infty$.  Hypothesis RA is also relaxed in Section 10.

The case in which (\ref{eq3.6}) and (\ref{eq3.9}) hold is dubbed \emph{Hypothesis RS} (short for risk-seeking).  In it, the assumption that $q^k$ is transient rules out bandits whose repeated play would earn expected utility of $+\infty$.

Hypothesis C supports nearly all of the results in this paper.
An exception occurs in Sections 7-8, where Hypothesis RN (and only it) is shown to accommodate constraints that link the bandits.

\subsection{Transient matrices}

A central role is played by matrices that are nonnegative and transient.  Relevant information about these matrices is contained in Proposition \ref{prop2.1-111210}, below.  It employs this nomenclature; vectors  $x$  and  $y$  that have the same number of entries satisfy   $x \gg y$  if and only if $x_j >  y_j$ for each $j$.

\begin{prop} \label{prop2.1-111210}
Let  $Q$ be a nonnegative  $n\times n$   matrix.  The following are equivalent:
\bi
\item[(a)]	The matrix  $Q$  is transient.
\item[(b)]	The matrix  $(I - Q)$  is invertible, and $(I-Q)^{-1} = I+Q+Q^2 + \cdots$.
\item[(c)]	There exists an  $n \times n$   vector  $f  \gg  0$  such that the equation  $(I-Q) x = f$   has a solution  $x \gg 0$.
\item[(d)]	There exists an $n \times 1$   vector  $y  \gg 0$  such that  $y \gg Q y$.
\ei
\end{prop}

\Pf
Demonstration that  $(a) \Rightarrow  (b) \Rightarrow   (c) \Rightarrow   (d) \Rightarrow  (a)$  is routine and is omitted. \EndPf

Parts of the analysis that follows could be simplified in the linear-utility case because a substochastic matrix  $q^k$   is transient if and only if termination occurs with positive probability after at most     $|N_k|$ transitions.

\subsection{Inheritance}

Hypothesis C is a property of the individual bandits.  Its implications for the multi-armed bandit are investigated next.  Let us recall that $S$  denotes the set of all multi-states of the multi-armed bandit.  Each stationary nonrandomized policy $\pi$  has an $|S| \times |S|$ \emph{transition rate matrix} $Q^{\pi}$    that is given, for each pair  $s$  and  $t$  of multi-states, by
\beq Q^{\pi} (s,t) = \lcb \begin{array}{ll} q(s_{\pi(s)},j) \quad & \mbox{if } \ t=s_{\setminus \pi (s)} \cup \{j\} \\
0 & \mbox{otherwise.}
\end{array} \right . \label{eq3.10}
\eeq
Each stationary nonrandomized policy $\pi$  also has an $|S| \times 1$ \emph{reward vector} $R^{\pi}$
that is defined for each state  $s$  in  $S$  by
\beq R^{\pi} (s) = r(s_{\pi(s)}) \, . \label{eq3.11}
\eeq

\begin{prop} \label{prop2.2-111210}		  .
Consider any stationary nonrandomized policy  $\pi$.  Condition (\ref{eq3.6}) guarantees that  $Q^{\pi}$   is nonnegative and transient.
\end{prop}

\Pf
\emph{(adapted from \cite{DPR})}.  By hypothesis, each bandit  $k$  has a transition matrix $q^k$    that is nonnegative and transient.  That  $Q^{\pi}$   is nonnegative is immediate from (\ref{eq3.10}).  Part (d) of Proposition \ref{prop2.1-111210} guarantees that each bandit  $k$  has a column vector  $x^k \gg 0$  such that $x^k \gg q^k x^k$.  Denote as  $y$  the  $|S| \times 1$   vector whose entry $y_s$  for multi-state  $s$  is given by $y_s = x_{s_1}^1 x_{s_2}^2 \dots x_{s_k}^K$.  It is clear that  $y \gg 0$.  Consider any multi-state  $s$;  set  $k  = \pi(s)$  and set  $i  = s_k$.  The nonzero entries in the $s^{\rm th}$ row of $Q^{\pi}$   correspond to the nonzero entries in the  $i^{\rm th}$ row of $q^k$,  and the inequality  $x^k \gg q^k x^k$  guarantees $y_s > [Q^{\pi} y]_s$.  This holds for each multi-state $s$,  so part (d) of Proposition \ref{prop2.1-111210} guarantees that $Q^{\pi}$  is transient.  \EndPf

That (\ref{eq3.6}) is inherited by the multi-armed bandit is the gist of Proposition \ref{prop2.2-111210}.  That (\ref{eq3.7})--(\ref{eq3.9}) are inherited is evident from (\ref{eq3.10}) and (\ref{eq3.11}).  Thus, the multi-armed bandit inherits the hypothesis that is satisfied by the individual bandits.

\subsection{A sequential decision process}
 	
A well-developed theory of sequential decision processes (c.f., Denardo~\cite{EVD67} or Veinott~\cite{vein69}) can be applied to the model whose local utility function is given by (\ref{eq3.1}) with transition rates that satisfy (\ref{eq3.6}).  Proposition \ref{prop2.2-111210} shows that each stationary nonrandomized policy $\pi$  has a transition rate matrix  $Q^{\pi}$   that is nonnegative and transient, so Part (b) of Proposition \ref{prop2.1-111210} shows that $(I-Q^{\pi})$  is invertible.  With the $|S| \times 1$    vector  $V^{\pi}$   defined by
\beq
V^{\pi} = (I-Q^{\pi})^{-1} R^{\pi} \qquad \forall \, \pi \in\Delta \, . \label{eq3.12}
\eeq
Part (b) of Proposition \ref{prop2.1-111210} also justifies the interpretation of the $s^{\rm th}$ entry in $V^{\pi}$   as the expected utility for starting in state  $s$  and using stationary nonrandomized policy $\pi$ until termination occurs.  Premultiplying (\ref{eq3.12}) by $(I-Q^{\pi})$    produces the familiar \emph{policy evaluation equation},
\beq
V^{\pi} = R^{\pi} + Q^{\pi} V^{\pi}. \label{eq3.13}
\eeq		 .

With the $|S| \times 1$   vector $F$  defined by
\beq F(s) = \max \lcb V^{\delta} (s): \delta \in \Delta \rcb \qquad \forall \, s \in S \, ,
\label{eq3.14} \eeq
the number  $F(s)$  equals the largest expected utility obtainable from any stationary nonrandomized policy, given starting state $s$.  A policy $\pi$  is said to be \emph{optimal}\/ if $V^{\pi} = F$.
The restriction to stationary  nonrandomized policies is justified because
Hypothesis C has been shown to suffice for such a policy to be
optimal over the class of all history-remembering policies,
see \cite{EVD67} or \cite{vein69}. Further, such a policy can be found by linear programming, by policy improvement, or by successive approximation.  None of these methods is practical when the number $|S|$ of multi-states is large, however.


\section{Labeling and data revision}\label{s4-111114}

Let us recall that each bandit $k$  has a distinct set $N_k$  of states, that $N$  is the union of all states of all bandits, that  $0$  is a special state that is not in  $N$,  and that termination is modeled by transition to state  $0$.  A \emph{labeling}  $L$  is the assignment to each $j \in N \cup \{0\}$  of a label  $L(j)$  that is an integer between  $1$  and  $|N|+1$,  with  $L(0) = |N| +1$    and with no two states having the same label.  Thus, each labeling  $L$  assigns a distinct label to each state in $N$,  and it assigns the highest label to state  $0$.

A stationary nonrandomized policy  $\pi$  for the multi-armed bandit is called a \emph{priority rule}\/ if it is determined by a labeling  $L$  like so:
\beq \label{eq4.1}		
\pi(s)  =  \argmin \lcb L(s_k) : 1 \leq k \leq K \rcb \qquad \forall \, s \in S \, .
\eeq
The priority rule $\pi$  in (\ref{eq4.1}) is said to be \emph{keyed}\/ to the labeling $L$.  Given any multi-state  $s$,  this priority rule plays the bandit  $k$  whose current state $s_k$ has the lowest label.

\subsection{Revised rewards and transition rates}

The notation is now simplified somewhat.  For the remainder of this section, bandit $k$  has  $n$  states (rather than  $|N_k|$   states), and these states are numbered  $1$  through  $n$.  This bandit's transition rates form the $n \times n$  matrix $q^k$,  and its rewards form the $n \times 1$  vector $r^k$.

Consider the state  $i$  in bandit  $k$  that has
\beq
i = \arg\min \{L(j): j \in N_k\} \, . \label{eq4.2}
\eeq
Suppose a multi-state  $s$  is observed that includes state  $i\in N_k$  and for which the priority rule  $\pi$  has $\pi(s)  =  k$.  The priority rule  $\pi$  will continue to call for bandit  $k$  to be played until it experiences a transition to a state other than  $i$.  This motivates the replacement of each transition rate  $q(j, p)$  and each reward  $r(j)$  in bandit $k$ by $\bar{q} (j,p)$ and $\bar{r} (j)$,  where:
\begin{eqnarray}
\bar{q} (i,p) & = & q(i,p) / [1-q(i,i)]  \qquad \qquad \mbox{ if } p \neq i \, , \label{eq4.3} \\ 		
\bar{q} (j,p) & = & q(j,p) + q(j,i)\bar{q} (i,p)  \qquad \quad \mbox{ if } j \neq i \ \mbox{ and } p \neq i \, , \label{eq4.4} \\ 	
\bar{q} (j,i) & = & 0 \qquad \qquad \qquad \qquad \qquad \quad \ \mbox{ for each } j \, , \label{eq4.5} \\
\bar{r} (i) & = & r(i) / [1-q(i,i)] \, ,   \label{eq4.6} \\
\bar{r} (j) & = & r(j) + q(j,i) \bar{r} (i)  \qquad \qquad \quad \ \mbox{if } j \neq i \, .  \label{eq4.7}	 \end{eqnarray}
The selection of $i$ borrows from \cite{tsit94}, but that reference does not suggest any scheme
to update the data as is done in (\ref{eq4.3})-(\ref{eq4.7}).

Repeated play replaces the bandit's transition rate matrix  $q^k$ by the matrix $\bar{q}^k$ whose entries are given by (\ref{eq4.3})--(\ref{eq4.5}), and it replaces the bandit's reward vector $r^k$ by the vector $\bar{r}^k$ whose entries are given by (\ref{eq4.6})--(\ref{eq4.7}).  The revised transition matrix and reward vector are for a model in which transitions to state $i$  do not occur.  It will soon be demonstrated that $\bar{q}^k$ and $\bar{r}^k$  inherit the version of Hypothesis C that is satisfied by     $q^k$ and $r^k$.

\subsection{Elementary row operations}
 	
Equations (\ref{eq4.4}), (\ref{eq4.5}) and (\ref{eq4.7}) describe a model
in which the data of bandit $k$ has been revised so that no transitions occur to the state  $i$.
This process can be iterated.  The second execution of (\ref{eq4.2}) occurs with state  $i$  removed from $N_k$,  and it selects the state $\bar{i}$ in $N_k$ whose label is second lowest.  And so forth.  Algorithmically, the effect of repeated data revision is to begin with the $n \times (n+1)$  matrix (tableau) $[(I-q^k), r^k]$ and to use elementary row operations to alter the entries in this tableau like so:

\medskip

\noindent
{\bf Triangularizer} (for bandit  $k$  in accord with labeling  $L$).

\be
\item[1.]
Begin with the tableau $[(I-q^k), r^k]$.  Set  $M=N_k$.  While $M$ is nonempty, do Steps 2 and 3.
\item[2.]	Find the state  $i \in M$   whose label  $L(i)$  is smallest.  Set  $\alpha = 1/[1-q^k(i,i)]$.
\bi
\item[(a)]	Replace row $i$ of the tableau $[(I-q^k), r^k]$ by itself times the constant  $\alpha$.
\item[(b)]	For each state $j \in M \setminus \{i\}$, replace row $j$ of this tableau by itself plus the constant $q(j,i)$ times (the updated) row $i$; this update equates $q_{jt}$ to
    0 for $t=i$ and for each $t$ in $N\setminus M$.
\ei
\item[3.]	Replace  $M$  by $M \setminus \{i\}$.
\ee

The first execution of Step 2 replaces the tableau $[(I-q^k), r^k]$ by $[(I-\bar{q}^k), \bar{r}^k]$ where the entries in the $n \times 1$  vector  $\bar{r}^k$   and in the $n \times n$  matrix $\bar{q}^k$ are specified by (\ref{eq4.3})--(\ref{eq4.7}) with  $i$  as the state whose label is lowest.  The second execution of Step 2 replaces $\bar{q}^k$    by the transition rate matrix $\hat{q}^k$  for which transitions to state  $i$  are not observed and in which transitions to the state $\bar{i}$ whose label is second lowest are not observed, except for transition from  $i$  to $\bar{i}$.  And so forth.

\begin{prop} \label{prop3.1-111210}
Suppose that the data for bandit  $k$  satisfy Hypothesis RN, RA, or RS.  When the data for bandit  $k$  are triangulated in accord with a labeling  $L$,  each iteration of Step 2 produces a tableau $[(I-\bar{q}^k), \bar{r}^k]$  that satisfies the same hypothesis. \end{prop}

\Pf
By hypothesis, $q^k$ is nonnegative and transient.  The initial execution of Step 2 of the Triangularizer replaces $[(I-q^k), r^k]$ by  $[(I-\bar{q}^k), \bar{r}^k]$.  It does so by multiplying row $i$ by the positive number  $\alpha$  and then replacing each row $j$ other than $i$  by itself plus the nonnegative multiple $q(j,i)$   times the updated row $(i)$.  This guarantees $\bar{q}^k \geq 0$.  It further guarantees that $\bar{r}^k \leq 0$  if  $r^k \leq 0$  and that  $\bar{r}^k \geq 0$  if  $r^k \geq 0$.  In particular, (\ref{eq3.8}) and (\ref{eq3.9}) are preserved.
 	
Since $q^k$ is nonnegative and transient, Part (c) of Proposition \ref{prop2.1-111210} shows that there exists a vector  $f \gg 0$  such that the equation $(I-q^k) x = f$ has a solution  $x \gg 0$.  Let us apply the Triangularizer to the tableau $[(I-q^k),f]$.  The initial execution of Step 2 replaces  $[(I-q^k),f]$  by   $[(I-\bar{q}^k),\bar{f}]$.   Elementary row operations preserve the solutions to equation systems, so the strictly positive vector  $x$  satisfies $(I-\bar{q}^k) x = \bar{f}$.  Since $\bar{q}^k \geq 0$,  part (c) of Proposition \ref{prop2.1-111210} also shows that $\bar{q}^k$ is transient, hence that (\ref{eq3.6}) is preserved.
 	
Finally, suppose that $q^k$ satisfies (\ref{eq3.7}).  With  $e$  as the $n \times 1$  vector of  1's,  note that $(I-q^k)e=g$ with $g \geq 0$.  As noted above, $(I-\bar{q}^k)e = \bar{g}$ with  $\bar{g} \geq 0$,  which shows that (\ref{eq3.7}) is preserved.

It has been demonstrated that Hypotheses RN, RA and RS are preserved by the first execution of Step 2 of the Triangularizer.  Iterating this argument completes the proof.      \EndPf

The computational effort for executing the Triangularizer is determined in the next result.

\begin{prop} \label{prop3.2-111210}
With $n\equiv |N_k|$, executing the Triangularizer on bandit $k$ entails $\frac{2}{3}|n|^3-\frac{1}{2}|n|^2+\frac{2}{3}|n|$ arithmetic operations.
\end{prop}
\Pf
The computation of the $[1-q^k(i,i)]$'s is Step 1 requires  $n$ subtractions.
Next consider
the execution of Step 2 when $|M|=m$. As the entries in row $i$ indexed by the columns
of $N_k\setminus M$ are zero and are not changed in Substep 2(a) and as
$\frac{1-q^k(i,i)}{1-q^k(i,i)}=1$, Substep 2(a) requires  $m$ divisions (including
the update or $r^k(i)$).
Also, Substep 2(b) requires $(m-1)m$ additions and  $(m-1)m$ multiplications.
Thus, the total number of arithmetic operation needed
to execute both substeps is $m+2m(m-1)=2m^2-m$.
As $\sum_{m=1}^n m^2=\sum_{m=1}^n 2{m\choose 2} + {m \choose 1}=2{{n+1}\choose 3}+{{n+1}\choose 2}
=\frac{1}{3}n^3+\frac{1}{2}n^2+\frac{1}{6}n$
and $\sum_{m=1}^n m=\frac{n^2+n}{2}$,
the total number of arithmetic operations needed to execute the triangularizer is
$n+2[\frac{1}{3}n^3+\frac{1}{2}n^2+\frac{1}{6}n]-\frac{n^2+n}{2}
=\frac{2}{3}n^3 -\frac{1}{2}n^2+\frac{2}{3}n$.
 \EndPf

 \subsection{Illustration}
	
The net effect of the Triangularizer is easiest to visualize when state 1 has the lowest label, state 2 has the next lowest label, and so forth.  In this case, the Triangularizer transforms the tableau $[(I-q^k), r^k]$ into the $n \times (n+1)$  tableau  $[(I-\tlq^k), \tlr^k]$   whose entries have the format,
\beq
\lsb \begin{array}{cccccc}
1 & - \tlq(1,2) & - \tlq (1,3) & \cdots & - \tlq (1,n) & \tlr (1) \\
0 & 1 & - \tlq (2,3) & \cdots & - \tlq (2,n) & \tlr (2) \\
0 & 0 & 1 & \cdots & - \tlq (3,n) & \tlr (3) \\
\vdots & \vdots & \vdots & \ddots & \vdots & \vdots \\
0 & 0 & 0 & \cdots & 1 & \tlr (n) \end{array} \rsb \ , \label{eq4.11}
\eeq
with  1's  on the principal diagonal and  0's  below that diagonal.
With finalized data,
each transition is to a state having a larger label, and termination is
guaranteed to occur after  $n=|N_k|$  transitions.

With linear utility---but not with exponential utility---the finalized data have simple interpretations:  Given that bandit  $b(i)$  is in state  $i$,  the number $\tlr(i)$    equals the expectation of the income that will be earned if bandit  $k$  is played until it  experiences transition to a state whose label exceeds  $L(i)$,  and  $\tlq(i,j)$   is the probability that this transition will occur to state  $j$.

\subsection{Finalized data}

Here and henceforth, tildes are used to identify the rewards and transition rates with which the Triangularizer ends, as in $\tilde{r}(i)$, $\tlq (i,j)$, $\tlr^k$,  and $\tlq^k$,  and these data are said to be \emph{finalized}.  The data for state  $i$  reach their finalized values when Step 2 is executed for state  $i$.  In other words, after Step 2 is executed for state  $i$,  no further changes occur in the $i^{\rm th}$ row or column of the tableau $[(I-q^k),r^k]$.

Let $\pi$ be the priority rule that is keyed to the labeling  $L$.  Equations (\ref{eq3.10}) and (\ref{eq3.11}) specify the $|S| \times |S|$ matrix $Q^{\pi}$ and the $|S| \times 1$  vector $R^{\pi}$ in terms of the original data.  Their analogs $\tlQ^{\pi}$ and $\tlR^{\pi}$  using finalized data are:
\begin{eqnarray}
\tlQ^{\pi} (s,t) & = & \lcb \begin{array}{ll} \tlq (s_{\pi(s)}, j) \quad & \mbox{if } t = s_{\setminus \pi(s)} \cup \{j\} \\
0 & \mbox{otherwise} \end{array} \right . \ , \label{eq4.8} \\[0.2cm]
\tlR^{\pi} (s) & = & \tlr (s_{\pi(s)}) \, . \label{eq4.9}
\end{eqnarray}       	

It was demonstrated in Section 3 that Hypothesis C is inherited by the multi-armed bandit.  Hence, with $V^{\pi}(s)$ as the expected utility for starting in state  $s$  and using priority rule  $\pi$,  the vector $V^{\pi}$  is the unique solution to $V^{\pi} = R^{\pi} + Q^{\pi} V^{\pi}$.  Proposition \ref{prop3.1-111210} shows that the model with finalized data also inherits Hypothesis C, hence that its reward vector $\tilde{V}^{\pi}$   is the unique solution to the policy evaluation equation
\beq
\tlV^{\pi} = \tlR^{\pi} + \tlQ^{\pi} \tlV^{\pi} \, . \label{eq4.10} \eeq	
That finalizing the data preserves expected utility is the gist of:

\begin{prop}
\label{prop3.3-111208}
Suppose Hypothesis C is satisfied.  Let $\pi$  be a priority rule that is keyed to a labeling $L$.  Then    $\tlV^{\pi} = V^{\pi}$.
\end{prop}

\Pf
A sequence of elementary row operations akin to those in the Triangularizer transforms the tableau $[(I - Q^{\pi}), R^{\pi}]$ into $[(I - \tlQ^{\pi}), \tilde{R}^{\pi}]$.  Elementary row operations preserve the set of solutions to an equation system.  Hence, since $V^{\pi}$ is the unique solution to  $(I-Q^{\pi}) V^{\pi} = R^{\pi}$,  it is the unique solution to  $(I-\tlQ^{\pi})v^{\pi} = \tilde{R}^{\pi}$. \EndPf

The Triangularizer first appeared in \cite{DPR}, with an elaborate analysis.  An antecedent to it appeared in Kaspi and Mandelbaum~\cite{kasp98}, and a contemporaneous account can be found in Sonin~\cite{soni08}.  That elementary row operations simplify the analysis seems not to have been observed previously, however.

\section{Policy evaluation}\label{s5-111117}

Throughout this section, $\hat{s}$ is any given multistate,
 $L$ is any  given labeling and $\pi$  is the priority rule
 that is keyed to $L$.
An algorithm that computes the expected utility $V^{\pi} (\hat{s})$
will be presented. Proposition \ref{prop3.3-111208} shows that $V^{\pi}=\tlV^{\pi} $
for which reason $V^{\pi} (\hat{s})$ can - and will be - computed using
finalized data. With finalized data, for (any given) state $i$ bandit $b(i)$ is played
at most once while in state $i$ and the finalized return $\tlr (i)$ is earned if that event occurs. The expected utility $\tlV^{\pi}(\hat{s})$ is then a linear combination of the finalized rewards say
\beq\label{eq5.7}
V^{\pi} (\hat{s}) = \sum_i z(i) \tlr (i) \, ;
\eeq
in the case of linear utility, $z(i)$ is the probability that bandit $k$ is played when its state is $i$, with finalized data.

A recursion will be used to compute the $z(i)$'s. Each step of this recursion updates entries in a set of vectors, on per bandit. For $p=1,\dots,K$, the vector $y^p$ has $|N_p|$ entries, one per state, and is initialized by
\beq
y^p (j) = \lcb \begin{array}{ll}
1 \qquad & \mbox{ if } \ j = \hat{s}_p \\
0 & \mbox{ if } \ j \in N_p \setminus \{\hat{s}_p\} \end{array} \right . . \label{eq5.1}
\eeq
Successively, for  $n  =  1, 2, \dots, |N|$,  this procedure selects the state  $i$  having  $L(i) = n$,  sets  $k  =  b(i)$,  updates  $y^k$   by
\begin{eqnarray}
y^k (j) & \leftarrow & [y^k (j) + y^k (i) \tlq (i,j)] \qquad \mbox{ if } \ j \in N_k \setminus \{i\} \, , \label{eq5.2} \\[0.2cm]
y^k (i) & \leftarrow & 0 \, , \label{eq5.3}
\end{eqnarray}
and makes no change in  $y^p$   for any $p \neq k$.  Equation (\ref{eq5.2}) augments the transition rate $y^k (j)$   to state  $j$  by the transition rate $y^k (i) \tlq (i,j)$   to state  $i$  and then directly to  $j$.  Equation (\ref{eq5.3}) reflects the fact that no state  $i$  is revisited when finalized data are employed.

The analysis of this procedure is eased by defining, for  $n  =  1, 2, \dots, |N|$,
\beq
P_n = \{s \in S: n > \min \{L (s_k): 1 \leq k \leq K\} \, . \label{eq5.4}	
\eeq
Evidently,  $P_n$   contains those multi-states that include a state whose label is less than $n$.

\begin{prop}
\label{prop4.1-111210}
Suppose Hypothesis C is satisfied.  Interrupt the execution of (\ref{eq5.2})--(\ref{eq5.3}) just prior to the iteration in which it selects the state  $i$  having  $L(i)  =  n$.  At this moment, the quantity $y^p (j)$ equals the aggregate transition rate with finalized data of bandit  $p$  from state $\hat{s}_p$  to state  $j$  due to play at each multi-state in $P_n$.
\end{prop}

\Pf
When  $n  =  1$,  this result  corresponds to the initial conditions.  Suppose it holds for  $n\ge 1$.  Expressions (\ref{eq5.2}) and (\ref{eq5.3}) show that it holds for  $n+1$. \EndPf

\vspace{-.25cm}
Proposition \ref{prop4.1-111210} prepares for the analysis of the: 	

\medskip

\noindent
{\bf Evaluator} (for starting multi-state $\hat{s}$, labeling  $L$ and priority rule $\pi$ that is keyed to $L$).
\be
\item[1.]	For each bandit  $k$,  define  $y^k$  by (\ref{eq5.1}).  Set  $V  =  0$  and  $n  =  1$.  While  $n \leq |N|$,  do Steps 2 and 3.
\item[2.]	Let  $i$  be the state whose label  $L(i)$  equals  $n$,  and set  $k  =  b(i)$.  Replace  $V$  by
\beq V + \tlr (i) y^k (i) \prod_{p \neq k} \lsb \sum_{j \in N_p} y^p (j) \rsb \, . \label{eq5.5}
\eeq
\item[3.] Execute (\ref{eq5.2}) and then (\ref{eq5.3}) for bandit $k$.  Then replace $n$  by  $n+1$.
\ee

\medskip
The next result shows that the Evaluator determines $V^{\pi} (\hat{s})$.

\begin{prop} \label{prop4.2-111210}
Suppose Hypothesis C is satisfied.  The Evaluator terminates with $V  = V^{\pi} (\hat{s})$.
\end{prop}

\Pf
For $i\in N$, set  $n  =  L(i)$  and  $k  =  b(i)$.  The coefficient  $z(i)$  in (\ref{eq5.7}) equals the aggregate transition rate from multi-state $\hat{s}$ to the set of multi-states $s$  that have  $n = \min \{L(s_p): 1 \leq p \leq K\}$.  From Proposition \ref{prop4.1-111210}, we obtain
$ z(i)  = y^k (i) \prod_{p \neq k} \lsb \sum_{j \in N_p} y^p (j) \rsb $,
which completes the proof. \EndPf

\medskip
The computational effort for executing the Evaluator is determined in the next result.

\begin{prop} \label{prop4.3-111210}
With $n\equiv \sum_{k=1}^K |N_k|$, executing the Evaluator entails $\sum_{k=1}^K \frac{3}{2}|N_k|^2+\frac{7}{2}n-5$ arithmetic operations
 (beyond the effort required to apply the Triangularizer on each bandit).
\end{prop}
\Pf
Augment the evaluator by keeping a record of $w^p=$ \ $ \sum_{j\in N_p} y^p(j)$ for each $p$
and of $w=\prod_{p=1}^{|K|} \lsb \sum_{j \in N_p} y^p (j)\rsb$. The initial value
of each of these expressions is 1. Keeping record of these expression
will facilitate the
computation of the bracketed terms in (5.5) by a single division.

Consider the implementation of Step 2 when $i\in N_k$ is selected and
$m$ is the number of states in $N_k$ whose label is lager than $L(i)$. In
this case the execution of (5.5) is Step 2 requires one addition, 2
multiplications and one division, totalling 4 arithmetic operations.
Also, in step 3, (5.2) has to be implemented only to the $m$ states in
$N_k$ whose label is higher than
$L(i)$, requiring $m$ additions and $m$ multiplications,
totalling $2m$ arithmetic operations.
Next, $\sum_{j\in N_p} y^p(j)$ has to be updated only for $p=k$
and this update requires
$m-1$ additions. Also, the update of
$\prod_{p=1}^{|K|} \lsb \sum_{j \in N_p} y^p (j)\rsb$
requires the multiplication of the old value by the ratio of the new and old values of
$\sum_{j\in N_k} y^k(j)$, requiring 2 arithmetic operations.
The total number of arithmetic operation
applied to execute steps 2 and  3 over all states $i$ is then
\[
\sum_{k=1}^K \sum_{m=1}^{|N_k|-1}(3m+5)=\sum_{k=1}^K \frac{(3|N_k|+10)(|N_k|-1)}{2}
=\sum_{k=1}^K \frac{3}{2}|N_k|^2+\frac{7}{2}n-5 .\hspace{2cm} \blot
\]

\medskip
To our knowledge, the computation of $V^{\pi}(\hat{s})$  for
a particular priority policy $\pi$ and particular starting state $\hat{s}$ is new.  With a different function (\ref{eq5.1}), the Evaluator and its work bound apply to any initial distribution over the multi-states that is in product form (except that the initial values of the $w^p$'s and $w$ of the proof of Theorem
5.2 will require n-1 additional arithmetic operations).

\section{Pairwise comparison and preference}\label{s6-111117}

In this section, pairwise comparison will be used to identify a state that is ``best" amongst a group of states, and the data for that state's bandit will be revised accordingly.  The \emph{amplification}  $a(i)$  of state  $i$  is now defined by
\beq
a(i) = \sum_{j \in N_{b(i)}} q(i,j) \, . \label{eq6.1}
\eeq			
Under Hypothesis RN, each amplification is  $1$  or less.  In the risk-averse and risk-seeking cases, some states can have amplifications that exceed  $1$,  however.

Playing chain  $b(i)$  when its state is  $i$  earns reward  $r(i)$  and multiplies future rewards by the factor  $a(i)$.  State  $i $ is now said to be \emph{preferable}\/ to state  $j$  if
\beq	r(i) +  a(i) r(j)  >  r(j) + a(j) r(i) \, . \label{eq6.2}
\eeq
Suppose that state  $i$  is preferable to state  $j$:  if a multi-state  $s$  is observed that includes states  $i$  and  $j$,  playing bandit  $b(i)$  first and  $b(j)$  second is better than the other way around.  The definition of preference is applied even when states  $i$  and  $j$  are in the same bandit, however.

It will soon be seen that preference is not transitive, but that it can be refined in a way that is transitive.  To this end, states will be grouped into ``categories."  The rule by which a category is assigned to each state varies with the hypothesis.

\subsection{Categories under Hypothesis RN}

Under Hypothesis RN, each state  $j$  has  $a(j) \leq 1$,  and each state is assigned a category by this rule:
\bi
\item	\emph{Category 1}\/ consists of each state  $j$  that has  $a(j)  =  1 $ and  $r(j) \geq  0$.
\item \emph{Category 2}\/ consists of each state  $j$  that has  $a(j)  <  1$.
\item \emph{Category 3}\/ consists of each state  $j$  that has  $a(j)  =  1$  and  $r(j)  <  0$.
\ei
It is easy to see that each state  $j$  in category 1 that has  $r(j)  >  0$  is preferable to every state in category 2 and that each state in category 2 is preferable to every state in category 3.  But no state in category 1 is preferable to any state in category 3.  For this reason, preference is not transitive.  	
State  $i$  is now said to be \emph{weakly preferable}\/ to state  $j$  if the inequality,
\beq
r(i) +  a(i) r(j) \geq r(j) + a(j) r(i) \, , \label{eq6.3}
\eeq
holds strictly or if this inequality holds as an equation and the category of  $i$  is at least as small as the category of  $j$.  Under Hypothesis RN, each state  $i$  is assigned a \emph{ratio}\/  $\rho(i)$  by the following rule:
\beq
\rho (i) = \lcb \begin{array}{ll} + \infty & \mbox{ if state $i$ is in category 1,} \\
r(i) / [1-a(i)] \qquad \quad & \mbox{ if state $i$ is in category 2,} \\
- \infty & \mbox{ if state $i$ is in category 3.} \end{array} \right . \label{eq6.4}
\eeq
	 		
It is easy to check that state  $i$  is weakly preferable to state  $j$  if and only if $ \rho(i)  \geq \rho (j)$.  Evidently, weak preference is transitive.  A state $i$  that is weakly preferable to all others can be found with  $|N|-1$ comparisons.

\subsection{Categories under Hypothesis RA}

Under Hypothesis RA, a state  $j$  can have  $a(j) > 1$,  but each state  $j$  has  $r(j) \leq  0$,  and the states group themselves into categories like so:
\bi
\item	\emph{Category 1}\/ consists of each state  $j$  that has  $r(j)  =  0$  and  $a(j) \leq  1$.
\item \emph{Category 2}\/ consists of each state  $j$  that has  $r(j)  <  0$.
\item	\emph{Category 3}\/ consists of each state $j$  that has  $r(j)  =  0$  and  $a(i)  >  1$.
\ei
As before, state  $i$  is said to be \emph{weakly preferable}\/ to state  $j$  if  (\ref{eq6.3}) holds strictly or if (\ref{eq6.3}) holds as an equation and the category of state  $i$  is at least as small as the category of state  $j$.  Under Hypothesis RA, each state  $i$  is assigned a \emph{ratio}\/ $\rho(i)$  by this rule:
\beq
\rho (i) = \lcb \begin{array}{ll} + \infty \quad \quad & \mbox{ if state $i$ is in category 1,} \\
\lsb 1-a(i)\rsb  / r(i) \quad \quad & \mbox{ if state $i$ is in category 2,} \\
- \infty \quad \quad & \mbox{ if state $i$ is in category 3.} \end{array} \right . \label{eq6.5}
\eeq	
It is easy to check that  $i$  is weakly preferable to state  $j$  if and only if  $\rho(i) \geq \rho(j)$.

\subsection{Categories under Hypothesis RS} 	

In the risk-seeking case, each state  $j$  has  $r(j) \geq 0$,  and the states group themselves into categories by this rule:
\bi
\item \emph{Category 1}\/ consists of each state  $j$  that has  $r(i)  =  0$  and  $a(i) \geq 1$.
\item \emph{Category 2}\/  consists of each state  $j$  that has  $r(j)  >  0$.
\item \emph{Category 3}\/ consists of each state  $j$  that has $ r(j)  =  0$  and  $a(j)  <  1$.
\ei
State  $i$'s  \emph{ratio}\/ is now defined by:
\beq
\rho (i) = \lcb \begin{array}{ll} + \infty \quad \quad & \mbox{ if state $i$ is in category 1,} \\
\lsb a(i)-1 \rsb / r(i) \quad \quad & \mbox{ if state $i$ is in category 2,} \\
- \infty \quad \quad & \mbox{ if state $i$ is in category 3.} \end{array} \right . \label{eq6.6}
\eeq	
With this categorization, the definition of weak preference does not change.  Again, state  $i$  is weakly preferred to state  $j$  if and only if  $\rho(i) \geq \rho(j)$.

\subsection{Finding a weakly preferred state in a set}
The characterization of
``weakly preferred" under RN, RA and RS by comparing $\rho(\cdot)$
shows that the relation is transitive. Further, if
the $r(i)$'s and the $[1-a(i)]$'s for each state $i$ in a set
$U$ are available, then (\ref{eq6.4}), (\ref{eq6.5}) or (\ref{eq6.6}),
respectively, facilitate the
identification of a weakly preferred state in $U$ by applying at most
$|U|$ divisions and $|U|$ comparisons.

\subsection{A key result}

Proposition \ref{prop6.1} (below) would seem to have a simple proof, at least in the case of linear utility, but we are not aware of one.  The interested reader is referred to the proof of Theorem 5.2 in \cite{DPR}, which employs a delicate interchange argument.

\begin{prop} \label{prop6.1}
Suppose Hypothesis C is satisfied, and consider any state  $i$  that is weakly preferred to all others.  It is optimal to play bandit $ b(i)$  for every multi-state  $s$  that includes state  $i$.
\end{prop}

\subsection{Nomenclature}

In the discussion to follow, the data for bandits 1 through $K$ will be triangularized in parallel, rather than one after the other.  At any stage in that computation:
\bi
\item $r(j)$  and  $q(j, p)$  denote the \emph{current}\/ values of the data for state  $j$,
\item $\bar{r} (j)$   and $\bar{q} (j,p)$ denote values of the data after they have been \emph{updated}\/ by the next execution of Step 2 of the Triangularizer,
\item $\tlr (j)$ and $\tlq (j,p)$ denote the finalized values of the data.
\ei
The amplification for state  $j$  is denoted  $a(j)$, $\bar{a} (j)$  and $\tilde{a}(j)$  when it is given in terms of current, updated and finalized data, respectively.  The same is true of the ratio, $\rho (j)$.

It is recalled the data for state  $i$  attain their finalized values when Step 2 is executed for state  $i$.   Proposition \ref{prop6.2} (below) indicates how each execution of Step 2 of the Triangularizer affects the ratios.

\begin{prop} \label{prop6.2}
Suppose Hypothesis C is satisfied.  With  $M$  as any nonempty subset of $N_k$,  suppose that state  $i$  in bandit  $k$  be weakly preferable to the other states  $M$  with current values of the data for bandit  $k$.  Executing Step 2 of the Triangularizer for state  $i$  has these effects:
\begin{eqnarray}
\rho(i) & = & \tilde{\rho} (i) \, , \label{eq6.7} \\
\rho (i) & \geq & \bar{\rho} (j) \geq \rho(j) \qquad \qquad \forall \ j \in M \setminus \{i\} \, . \label{eq6.8}
\end{eqnarray}
\end{prop}

Equation (\ref{eq6.7}) states that finalizing the data for state  $i$  preserves its ratio.  Expression (\ref{eq6.8}) states that updating the ratio for a state  $j$  other than  $i$  can improve its ratio, but not above that for state $i$.  These observations are insightful, but Proposition \ref{prop6.2} is not used in this paper, and its proof is omitted.  Proposition \ref{prop6.2} facilitates the use of finalized data for each bandit, thereby enabling parallel computation.

\section{Optimization}\label{s7-111117}

 A labeling  $L$  is said to be \emph{optimal}\/ if the priority rule $\pi$ that
 is keyed to  $L$  has $V^{\pi} = F$, i.e., $\pi$ maximizes the expected
 utility that can be obtained from each starting multi-state.
  Proposition \ref{prop6.1} lays the groundwork for a variety of algorithms that
 identify an optimal labeling.  The Optimizer, which appears below, triangularizes the bandits contemporaneously, rather than one after the other.  Its first execution of Step 3 identifies the state $i$ that is weakly preferable to all others with respect to the original data.  Its first execution of Steps 3(a) and 3(b) update the data for bandit  $b(i)$  in accord with repeated play while in state $i$   and then remove state $i$.  The Optimizer then repeats Step 3 with updated data.  This recursion stops as soon as all states in one bandit have been removed

 \subsubsection*{Optimizer}

 \be
 \item[1.]	Begin  $C$  equal to the empty set.  For each bandit $k$,  insert in  $C$  a state $i \in N_k$    that is weakly preferable to every other state  $j \in N_k$  with respect to the original data.  Set  $n  =  1$.  For each bandit  $k$,  set  $M_k = N_k$.

\item[2.]	Do Step 3 while  $M_k$   is nonempty for each  $k$.

\item[3.]	Find a state  $i \in C$   that is weakly preferable to all other states in  $C$  with respect to current data.  Set  $k  =  b(i)$  and set  $L(i)  =  n$.  Then replace  $n$  by  $n+1$.

\bi
\item[(a)]	Use Step 2 of the Triangularizer to finalize the data for state  $i$  and to update the data for each state $j \in M_k \setminus \{i\}$.
\item[(b)]	Remove state  $i$  from  $C$.  Remove state  $i$  from $M_k$.  If $M_k$  is nonempty, insert in  $C$  a state $j \in M_k$   that is weakly preferable to all other states in  $M_k$   with respect to updated data.
\ei
 \ee

The Optimizer stops as soon as all of the states in any bandit have been labeled, with  $n - 1$  as the highest of the labels.  The unlabeled states can be assigned the labels  $n$  through $|N|$ in any way.  It will not matter: no bandit whose state is labeled  $n$  or higher will ever be played because it cannot have the lowest label.

The Optimizer applies the Triangularizer with respect to a labeling that
is determined on line.
At each stage, the state that gets the next label is selected so that it is
weakly preferred to
all states that have not yet been labeled, i.e., the states in $\cup_{k=1}^K M_k$.

\begin{prop}
\label{prop6.1-111210}
Suppose Hypothesis C is satisfied.  The Optimizer constructs a labeling  L  that is optimal.
\end{prop}

\Pf
Let  $i$  be the state selected at the initial execution of Step 3.  Weak preference is transitive, so Proposition \ref{prop6.1} shows that it is optimal to play bandit  $b(i)$  at each multi-state that includes state  $i$.  Setting  $L(i)  =  1$  is optimal.

Step 3(a) equates to   $0 $  the transition probability  $\tlq(j,i)$  for each state   $j $  in   $b(i)$,  and, for each state   $j $  in bandit   $b(i)$,  it updates the reward  $\tlr(j)$,  the transition probabilities $\tlq(j,p)$   for each  $p \neq i$  to account for repeated play while in state  $i$.

Step 3(b) removes state  $i$  from bandit  $b(i)$.  What remains is a multi-armed bandit with one fewer state.  Proposition \ref{prop3.1-111210} implies that the same version of Hypothesis C is satisfied by the  bandit with one fewer state.  Since weak preference is transitive, the state  $i$  that is selected at the second iteration of Step 3 is weakly preferable to all others in the model with revised data and one fewer state.  Proposition \ref{prop6.1} can be applied a second time, and state  $i$  can be assigned the label  $L(j)  =  2$.  Iterating this argument completes the proof. \EndPf

\medskip
The computational effort for executing the Optimizer is determined in the next result.

\begin{prop}
\label{prop6.2-111210}
With $n\equiv \sum_{k=1}^K |N_k|$, the Optimizer can be executed with $\frac{1}{2}\sum_{k=1}^K |N_k|^2+\frac{n}{2}$ arithmetic operations and  $\frac{1}{2}\sum_{k=1}^K |N_k|^2+n(K-\frac{1}{2})+{K\choose 2}$ comparisons plus the effort required to execute
the Trinangularizer on each bandit).
\end{prop}
\Pf
Augment the optimizer by recording a ranking of the elements of $C$
in decreasing weakly preferable order and corresponding
ratios of those states in $C$  that are in category 2.
The initial ranking can be accomplished with  $K\choose 2$ comparisons whereas the
ratios of the states in category 2 that enter $C$ in Step 1 are computed when the those
states are selected to enter $C$.

When state $i\in N_k$ gets a label, $M_k$ changes and the Triangularizer updates the
data of its states, including the $r(j)$'s and $[1-a(j)]$'s. At each stage,
finding a weakly preferred state in
$M_k$ can be accomplished with at most $|M_k|$ divisions (determining ratios for states
in category 2) and at most
$|M_k|-1$ comparisons. Updating the ranked list $C$ replaces the old state from $b(i)$ by $i$, requires at most $K-1$ comparisons.
So, the effort for executing the
Optimizer, beyond the effort required to execute
the Trinangularizer on each bandit, is bounded by $\sum_{k=1}^K\sum_{m=1}^{|N_k|} m=\frac{1}{2}\sum_{k=1}^K |N_k|^2+\frac{n}{2}$
arithmetic operations and $\sum_{k=1}^K\sum_{m=1}^{|N_k|} (K-1+m-1)
+{K\choose 2}=\frac{1}{2}\sum_{k=1}^K |N_k|^2+n(K-\frac{1}{2})+{K\choose 2}$ comparisons.
\EndPf

\medskip
Proposition \ref{prop6.1-111210}, \ref{prop3.1-111210}, \ref{prop4.3-111210} and \ref{prop6.2-111210}
show that an optimal priority rule and its expected utility  $F(s)$  for a particular starting state  $s$  can be computed with $\frac{3}{2}\sum_k |N_k|^3+O(N^2)=O(N^3)$ arithmetic operations and $O(N^2)$ comparisons. These last two bounds match the best existing bound for computing Gittins indices (obtained in \cite{NM07}, see \cite[p.43]{GGW11}).


This section is closed with the mention of an alternative to the Optimizer.  This alternative has two steps:  First, optimize within each individual bandit.  Second, use finalized data for each bandit and pair-wise comparison to rank the states $1$ through $|N|$ by weak preference.  Proposition \ref{prop6.2} shows that the priority rule that is keyed to this ranking (labeling) is optimal.  This procedure also requires work proportional to $\sum_k |N_k|^3$.

\section{Optimization with Constraints}\label{s8-111117}

For the case of a linear utility function that satisfies Hypothesis RN, the multi-armed bandit is now generalized to include a finite number $W$  of constraints, each on a particular type of reward.  Including the objective, there are now  $W+1$  types of reward, which are numbered  $0$  through  $W$.  The objective measures \emph{type-0}\/ reward and the $w^{\rm th}$ constraint places a lower bound $C_w$ on the expected \emph{type-w}\/ reward.

The initial multi-state  $s$  is given, and the object is to maximize the expectation of the type-0 reward subject to constraints that, for each  $w$,  keep the expectation of the type-$w$ reward is at least as large as $C_w$.  The main thrust of this section is to use column generation to construct an optimal solution to the constrained problem that is an initial randomization over  $W+1$  priority rules.  At the end of the section, the approach taken here is compared with a more classic one.

It is known (c.f., Feinberg and Rothblum~\cite{FR11}) that an optimal policy can be found among the initial randomizations over stationary deterministic policies.  This lets the multi-armed bandit problem with constraints be formulated as:

\medskip

\noindent
\begin{tabbing}
{\bf Program 1}.  Maximize  \= $\sum_{\delta} \alpha^{\delta}V_0^{\delta} (s)$, \= subject to the \= constraints \\[0.2cm]
\> $\sum_{\delta} \alpha^{\delta}$\> $=1$, \\[0.2cm]
\> $\sum_{\delta} \alpha^{\delta} V_w^{\delta} (s)$ \> $\geq C_w$ \> for \ \ $w = 1,2,\dots, W$,\\[0.2cm]
\> \ \quad $\alpha^{\delta} \geq 0$ \>\> for all  \ $\delta$,
\end{tabbing} 	
where it is understood that the sum is taken over all stationary deterministic policies  $\delta$  and where  $V_w^{\delta}(s)$   denotes the expectation of the type-$w$ utility that is earned if one starts at
multi-state  $s$  and uses policy  $\delta$.  Program 1 has only  $W+1$  constraints, but it can have a gigantic number of decision variables (one for each stationary deterministic policy $\delta$  and one for each slack variable), and its data include the type-$w$ reward  $V_w^{\delta} (s)$  for each  $w$  and each policy  $\delta$.

Program 2, below, is in the same format as Program 1.  Program 2 has one decision variable for each priority rule  $\pi$,  rather than for each policy $\delta$.

\medskip

\noindent
\begin{tabbing}
{\bf Program 2}.  Maximize  \= $\sum_{\pi} \alpha^{\pi}V_0^{\pi} (s)$, \= subject to the \= constraints \\[0.2cm]
$ \ \qquad y_0:$ \> $\sum_{\pi} \alpha^{\pi}$\> $=1$, \\[0.2cm]
$ \ \quad -y_w:$ \> $\sum_{\pi} \alpha^{\pi} V_w^{\pi} (s)$ \> $\geq C_w$ \> for \ \ $w = 1,\dots, W$,\\[0.2cm]
\> \ \quad $\alpha^{\pi} \geq 0$ \>\> for all  \ $\pi$.
\end{tabbing} 	

There are fewer priority rules than polices, but the number of priority rules can still be enormous.  Multipliers have been assigned to the constraints of Program 2.  These multipliers will be used in column generation.

\subsection{Preview}
 	Although Program 2 has fewer columns than does Program 1,
 computing the data it requires would still be onerous.
 Much of this computation can be avoided by coupling the
 simplex method with column generation.
 To indicate how, we suppose that a feasible basis for Program 2 has been found.
 This feasible basis consists of $W+1$ columns (the constraint matrix has full rank).
 It prescribes value of the basic variables and of the multipliers $y_0$ and $-y_1,\dots,-y_W$.  These multipliers are used to define rewards in an
 unconstrained bandit problem whose optimal solution
 (found by the Optimizer) identifies a
 priority rule $\lambda$  whose corresponding column has reduced cost (marginal profit) $c^\lambda$ that is the largest.  If  $c^\lambda$ equals zero, the current basis is optimal.  Alternatively, if  $c^\lambda$   is positive, the Evaluator is used to compute the coefficients  $V_0^\lambda,\dots,V_W^\lambda$.  A simplex pivot is then executed, and the process is repeated.

\subsection{Feasibility}

Each column of Program 2 is a column of Program 1.  Thus, if Program 2 is feasible, Program 1 must also be feasible.  The converse is established in:

\begin{prop} \label{prop7.1-111210}
Suppose Hypothesis RN is satisfied.  If Program 1 is feasible, Program 2 is also feasible.
\end{prop}

\Pf
We will prove the contrapositive.  Suppose that Program 2 is not feasible.  An application of Farkas' lemma (equivalently of the duality theorem of linear programming) shows that there exist numbers  $y_0$  and  $y_1$ through  $y_w$ such that
\begin{eqnarray}
y_0 - \sum_{w=1}^W y_w V_w^{\pi} (s) &  \geq & 0 \qquad \mbox{ for all } \ \pi \, , \label{eq8.1}  \\
y_w & \geq & 0 \qquad \mbox{ for } \ w = 1,\dots, W \, , \label{eq8.2} \\
y_0 - \sum_{w=1}^W y_w C_w & < & 0 \, . \label{eq8.3}
\end{eqnarray}

The numbers $y_1$  through  $y_W$  will be used as weights for the rewards $r_1(i)$  through $r_W(i)$.   Consider an unconstrained multi-armed bandit in which the reward  $R(i)$  for playing bandit  $b(i)$  while its state is $i$  is given by  $R(i)=y_1 r_1(i) + \cdots + y_W r_W (i)$.  Expression (\ref{eq8.1}) states that with reward  $R(i)$  for each state  $i$,  no priority rule $\pi$  has aggregate reward that exceeds $y_0$.  Proposition \ref{prop6.1} shows that a priority rule is optimal.  Thus,
\beq
y_0 - \sum_{w=1}^W y_w V_w^{\delta} (s) \geq 0 \qquad \quad \mbox{for all } \ \delta \, . \label{eq8.4}
\eeq
where  $\delta$ ranges over all stationary deterministic policies.  A solution exists to (\ref{eq8.2})--(\ref{eq8.4}), so a second application of Farkas' lemma shows that no solution exists to the constraints of Program 1.  \EndPf

Thus, Program 1 is feasible if and only if Program 2 is feasible.  Phase I of the simplex method will be soon used to determine whether Program 2 is feasible and, if so, to construct a feasible basis with which to initiate Phase II of the simplex method.  For the moment, it is assumed that a feasible basis for Program 2 has been found.

\subsection{Phase II}

The constraint matrix for Program 2 includes a column for each of the  $W$  slack variables.  These columns are linearly independent of each other, and they are linearly independent of the other columns.  Thus, the rank of its constraint matrix equals the number $W+1$  of its rows, and each basis for Program 2 consists of exactly  $W+1$ columns.  Let us consider an iteration of Phase II.  At hand at the start of this iteration is a feasible basis, its basic solution and its multipliers.  This information includes:
\bi
\item The data (column) for each of the  $W+1$  basic variables.
\item The basic solution (the $\alpha^{\pi}$'s) for this basis.
\item The multipliers  $y_0$   and  $y_1$    through $y_W$    for this basis.
\ei
The multipliers  $y_1$   through $y_W$   are nonnegative, and each priority rule  $\lambda$  has reduced cost $\bar{c}^{\lambda}$  that is given by
\[
\bar{c}^{\lm} = V_0^{\lm} (s) + \sum_{w=1}^W y_w V_w^{\lm} (s) - y_0 \, . \]

Computation of the reduced cost of each nonbasic priority rule $\lm$ would be an onerous task, but it is not necessary.  To determine whether or not the current basis is optimal and, if not, to find a priority rule has the largest (most positive) reduced cost, one can solve the \emph{unconstrained}\/ multi-armed bandit problem with the reward  $R(i)$   for playing bandit  $b(i)$  when its state is  $i$  given by
\beq R(i) = r_0 (i) + \sum_{w=1}^W y_w r_w (i) \, . \label{eq8.5}
\eeq
	
With these rewards, the Optimizer in Section 7 computes a priority rule  $\pi$  that is optimal.  Also, the Evaluator in Section 5 computes the expected return $V^{\pi} (s)$ for starting in multi-state $s$  and using this priority rule.  If $V^{\pi} (s) \leq y_0$,  no nonbasic variable has a positive reduced cost, so the current basis is optimal.  If $V^{\pi} (s) > y_0$,  the column for priority rule $\pi$  enters the basis.  To compute $V_w^{\pi} (s)$ for  each $w$, use the Triangularizer and Evaluator for priority rule  $\pi$.  In this computation, the finalized rewards vary with $w$ but the $y^k(j)$'s do not.  To complete an iteration of Phase II, execute a feasible pivot with $\alpha^{\pi}$  as the entering variable.

\subsection{Phase II recap}
	
Each feasible basis and its basic solution prescribe an initial randomization (with weight $\alpha^{\pi}$    assigned to priority rule $\pi$)  over $(W+1-p)$    priority rules, where  $p$  equals the number of slack variables that are basic.

The multipliers for the current basis determine the data of an unconstrained bandit problem, and the procedure in prior sections computes its optimal priority rule $\pi$    and its expected return, $V^{\pi}(s)$.  If  $V^{\pi}(s)$   does not exceed $y_0$,  the current basis is optimal.
If  $V^{\pi}(s)$   exceeds $y_0$,  the Evaluator is used to compute the coefficients  $V_0^{\pi}$ through     $V_W^{\pi}$ of the entering variable.  A simplex pivot is then executed.

The pivot itself requires work proportional to $(W+1)^3$.  Identifying the entering variable and its column of coefficients entails work proportional to $(W+2) [\sum_k|N_k|^3]$.  Only a few iterations may be needed to find a good basis, or an optimal basis, but that is not guaranteed.

\subsection{Phase I}

It remains to determine whether or not Program 2 is feasible and, if it is feasible, to construct a feasible basis with which to initiate Phase II.   These tasks will be accomplished by ``bringing in" the constraints of Program 2, one at a time.  Starting with  $n  =  1$,  the $n^{\rm th}$ iteration of Phase I is initialized with a randomization over  $n-1$  priority rules that satisfy the
first  $n-1$ constraints.  The $n^{\rm th}$ iteration maximizes type-$n$ reward, using the Phase II column generation scheme described above.  If the type-$n$ income can be made as large as $C_n$, a basis has been found with which to initiate the $n+1^{\rm st}$ iteration.  If not, no feasible solution exists to Program 2.    	

\subsection{The classic formulation}

An optimal policy for a discounted Markov decision problem with  $W$  constraints
can be found among the stationary randomized policies (c.f., Altman~\cite[page 102]{altm99}).
This can be accomplished by a linear program whose constraint matrix has one column per state-action pair, one row per state, and one row per constraint.  The multi-armed bandit has  $J = \prod_{k=1}^K |N_k|$
   multi-states  and  $K$  actions per multi-state.  Its constraint matrix has  $W+J$   rows and $K \times J$ columns.  The classic formulation has fewer columns than does Program 2, but it has many more rows.
 	
The classic formulation can also be attacked by column generation, but doing so would be unattractive because the formulation would have more columns and many more rows than the one we propose.

\subsection{A roadblock}

With a linear utility function, multiple types of reward can be handled by column generation and by the classic method.  Both methods utilize the fact that each transition rate  $q(i, j)$  is independent of the reward type.
 	
Let us consider what occurs when multiple types of rewards are introduced in the model with exponential utility.  Note from (\ref{eq3.4}) and (\ref{eq3.5}) that the payoff $x(i,j)$ appears in the formula for the transition rate $q(i,j)$.  Having multiple types of income causes the transition rate $q(i,j)$ to vary with the reward type.  Consequently, our column generation method (and the classical one) can be applied only
 when the type-$w$ payoff $x_w (i,j)$ is independent of  $w$, for instance, this is the case when income is earned only at termination.

\section{Structural properties}\label{s9-111117}

In the prior section, it was demonstrated that an optimal policy for a constrained bandit problem can be found among the initial randomization over  $W+1$  priority rules.  In the current section, the structure of this optimal policy is probed.

A transient Markov decision problem (MDP) with  $W$  constraints has an optimal solution that is an initial randomization over  $W+1$  deterministic policies $\delta^1$  through $\delta^{W+1}$   each of which differs from the next at precisely one state of the MDP; see Feinberg and Rothblum~\cite{FR11}.  When this MDP is a multi-armed bandit, these deterministic policies need not be priority rules, however.

Two labelings are now said to be \emph{adjacent}\/ if they are identical except that they exchange the states having labels  $k$  and  $k+1$  for exactly one value of  $k$.  The aforementioned property raises the question:  Does Program 2 have an optimal solution that is an initial randomization over priority rules that are keyed to a sequence of  $W+1$  labelings with the property that each labeling is adjacent to the next?  This question will be answered in the affirmative in the case of one constraint and in the negative in the case of more than one constraint.

\subsection{Adjacency with one constraint}

Let us consider a multi-armed bandit with one constraint.  We have seen that an optimal basis for Program 2 prescribes a randomization over at most two priority rules.  If its basic solution for this basis sets     $\alpha^j =  1$  for any  $j$,  only priority rule is used, and adjacency is trivial.

Let us denote as  $V_0^p$   and $V_1^p$    as the type-0 and type-1 utility for column $p$.  The case that requires analysis is that in which the optimal basis for Program 2 consists of columns  $j$  and  $k$  whose priority rules are keyed to different labelings.  For this to occur, the slack variable for the inequality constraint in Program 2 must not be basic, so this optimal basis assigns the columns  $j$  and  $k$  nonnegative values  $\alpha_j$   and  $\alpha_k$   that satisfy
\beq
\alpha_j V_1^j +\alpha_k V_1^k = C_1 \quad \mbox{ and } \quad \alpha_j + \alpha_k = 1 \, . \label{eq9.1}
\eeq
The optimal basis for Program 2 assigns to its constraints values of the multipliers  $y_0$   and $-y_1$  for which columns  $j$  and  $k$  have  $0$  as their reduced costs.  In other words,
\beq
0 = V_0^j + y_1 V_1^j - y_0 \quad \mbox{ and } \quad 0 = V_0^k + y_1 V_1^k - y_0 \, . \label{eq9.2}
\eeq
If  $V_1^j = V_1^k$,  equation (\ref{eq9.1}) guarantees that both columns have $C_1$ as their type-1 utility, and equation (\ref{eq9.2}) shows that both columns have the same type-0 utility, in which case it is optimal to play either column with probability 1, and a deterministic priority rule is optimal.

It remains to analyze the case in which  $C_1$   lies strictly between  $V_1^j$   and $V_1^k$.  The labelings to which columns  $j$  and  $k$  are keyed need not be adjacent, but columns  $j$  and  $k$  can be used to construct an optimal basis with labelings that are adjacent.  To indicate how, we
turn to the example in Table 1.  In this example, columns  $j$  and  $k$  assign identical labels to the states, except for the sets  $\{6, 7, 8, 9\}$  and  $\{13, 14\}$  of labels.

\bc
\begin{tabular}{ccccccccc}
\multicolumn{9}{c}{{\bf Table 1}.  An optimal basis.} \\
label	& $\dots$ & 6	& 7	& 8	& 9	& $\dots$	& 13	& 14 \\
\hline
column  $j$ &	$\dots$ &	a & b & c &	d	& $\dots$ & 	f	& g \\
column  $k$ &	$\dots$ &	d	& c	& b	& a	& $\dots$	& g	& f
\end{tabular}
\ec

Optimal solutions to the unconstrained multi-armed bandit having   $R(i) = r_0(i) + y_1 r_1 (i)$ for each state  $i$  are in product form.  As a consequence, every column  $p$  whose labeling permutes the labels assigned to the sets  $\{a, b, c, d\}$  and  $\{f, g\}$  of states has  $0$  as its reduced cost in Program 2.  A total of  $7  =  1 + 3 + 2 + 1$  interchanges of states whose labels are adjacent converts the permutation for column  $k$  into the permutation for column  $j$.  One of these interchanges must move the type-1 reward from the side of  $C_1$   on which $V_1^k$  lies to the side on which $V_1^j$   lies, and that switch identifies a pair of adjacent labelings.  This switch identifies a pair of priority rules that are keyed to adjacent labelings and whose columns form an optimal basis.  The pattern exhibited by this example holds in general.  The Triangularizer and Evaluator can be used to compute the reward vector for each labeling.

\subsection{Non-adjacency with two constraints}
	
For a multi-armed bandit problem with two constraints, an initial randomization over 3 priority rules has been shown to be optimal.  Examples exist in which no optimal solution is an initial randomization over priority rules that are keyed to a sequence of three adjacent labelings.  Such an example is now presented.  This example has 3 chains (bandits), each of which consists of a single state.  The three bandit's states are  $a$,  $b$  and  $c$,  respectively.  The multi-state  $(a, b, c)$  is observed initially.  Playing any bandit causes immediate termination.  Playing the bandit whose state is  $a$  earns the reward vector  $(1, 0, 0)$  whose entries are, respectively, the type-0, type-1 and type-2 reward.  Similarly, playing the bandit whose state is  $b$  earns reward vector  $(0, 1, 0)$,  and playing the bandit whose state is  $c$  earns reward vector  $(0, 0, 1)$.  The lower bounds on expected type-1 and type-2 rewards are $C_1=  0.3$  and $C_2 =  0.1$.  There are six labelings, which are listed below.  Labeling (v) has  $L(a) = 3$,  $L(b)  =  1$  and  $L(c)  =  2$,  for instance.

\bc
\begin{tabular}{lccc}
\ \quad state	& a	& b	& c \\
\hline
labeling (i) &	1 &	2 &	3 \\
labeling (ii) & 1 & 3 &	2 \\
labeling (iii)	& 2	& 1	& 3 \\
labeling (iv) & 2	& 3	& 1 \\
labeling (v) & 3 &	1	& 2 \\
labeling (vi) & 3	& 2	& 1
\end{tabular}
\ec

For this example, it is optimal to use labeling (i) or (ii) with probability of  $0.6$,  to use labeling (iii) or (iv) with probability $0.3$ and to use labeling (v) or (vi) with probability of  $0.1$.  But no sequence of three labelings, one from each pair, is adjacent.  For instance, labelings (i) and (iii) are adjacent to each other, but neither is adjacent to labeling (v) or (vi).   	

\section{Relaxing Hypothesis C} \label{s10-111117}

The model with a risk-averse exponential utility function can be generalized by replacing Hypothesis RA with these conditions:
\bi
\item Each bandit  $k$  has a transition rate matrix $q^k$ that is nonnegative.
\item At least one bandit  $k$  has a transition rate matrix $q^k$ that is transient.
\item	Every closed communicating class  $C$  of states in any bandit  $n$  has spectral radius of $(q^n)_{CC}$   that exceeds  $1$.
\ei

When Hypothesis RA is weakened in this way, the analysis becomes more intricate.  One difficulty stems from the fact that if a policy $\pi$  has a transition rate matrix  $Q^{\pi}$   that is not transient, its utility vector $V^{\pi}$  cannot satisfy (\ref{eq3.13}).  The fact that the risk-averse exponential utility function has  $u(0)  =  -1$  and the weakened hypothesis can be used to work around this difficulty by ruling out any stationary policy that plays a bandit at each state in any closed communicating class.  A second difficulty arises from the fact that the interchange argument in Proposition \ref{prop6.1}
can no longer rest on the classic results in \cite{EVD67} or \cite{vein69}.  The interested reader is referred to the analysis in \cite{DPR} and to the characterization of optimal policies in \cite{dena06}.
	
The linear-utility model can be generalized in a similar way.  It suffices that each bandit has a matrix $q^k$ that is substochastic, that at least one bandit  $k$  has a matrix $q^k$ that is transient, and that every closed communicating class of states in any bandit has a gain rate that is negative.

With each of these generalizations, only a minor change is required in the computation.  The change is to avoid playing bandit  $b(j)$  if at some point in the computation it has transition rate   $q(j, j)$  that equals or exceeds  $1$.

\section{Acknowledgements}\label{s11-111117}

The authors are pleased to acknowledge that this paper has benefited immensely from the reactions of Dr. Pelin Cambolat to earlier drafts.  The contribution of the second author has been supported in part by NSF grant CMMI-0928490. The contribution of the third author has been supported in part by ISF Israel Science Foundation) grant 901/10.

\end{document}